\documentclass[letterpaper]{article}
\usepackage{proceed2e}
\usepackage[margin=1in]{geometry}

\usepackage{algorithmic}
\usepackage{algorithm}

\usepackage[round]{natbib}
\usepackage{hyperref}
\usepackage{amsmath}
\usepackage{amsfonts}
\usepackage{amssymb}

\renewcommand{\P}{\mathbb{P}}

\newcommand{\bs}{\boldsymbol}

\newtheorem{definition}{Definition}
\newtheorem{theorem}{Theorem}
\newtheorem{lemma}{Lemma}
\newtheorem{corollary}{Corollary}
\newtheorem{proposition}{Proposition}
\newtheorem{example}{Example}

\newcommand*{\QEDA}{\hfill\ensuremath{\square}}

\def\cE{{\cal E}}

\def\cG{{\cal G}}

\def\cN{{\cal N}}

\def\cV{{\cal V}}

\usepackage{times}

\title{Identifiability in Gaussian Graphical Models}

\author{} 

\begin{document}

\maketitle

\begin{abstract}

In high-dimensional graph learning problems, some topological properties of the graph, such as bounded node degree or tree structure, are typically assumed to hold so that the sample complexity of recovering the graph structure can be reduced. With bounded degree or separability assumptions, quantified by a measure $k$, a $p$-dimensional Gaussian graphical model (GGM) can be learnt with sample complexity $\Omega(k \log p)$. Our work in this paper aims to do away with these assumptions by introducing an algorithm that can identify whether a GGM indeed has these topological properties without any initial topological assumptions. We show that we can check whether a GGM has node degree bounded by $k$ with sample complexity $\Omega(k \log p)$. More generally, we introduce the notion of a strongly $K$-separable GGM, and show that our algorithm can decide whether a GGM is strongly $K$-separable or not, with  sample complexity $\Omega(k \log p)$.
We introduce the notion of a generalized feedback vertex set (FVS), an extension of the typical FVS, and show that we can use this identification technique to learn GGMs with generalized FVSs.

\end{abstract}

\section{Introduction}





Probabilistic graphical models \citep{Pea:B88, Lau:B96, Whi:B90} have been increasingly studied as a means to represent relationships in multivariate distributions. In particular, the Gaussian graphical model (GGM) is a popular model with applications to many areas such as  object recognition and tracking \citep{Sud:06}, protein sequencing \citep{DurEddKroMit:99}, gene networks \citep{MohChuHanWitLeeFaz:12},  computer vision \citep{Isa:03} and neuroimaging \citep{Rya:12, BelVarBla:16}. Gaussian graphical models consists of a multivariate Gaussian distribution $\cN(\mu,\bs{\Sigma}$ and a graph $\cG$, where $\cG$ is defined by the Gaussian distribution.  We define an edge between two variables if and only if they are not conditionally independent given the rest of the variables. In the Gaussian setting, the edges between variables correspond to the non-zero entries of the precision matrix $\bs{\Omega} = \bs{\Sigma}^{-1}$, where $\bs{\Sigma}$ is the covariance matrix of the multivariate Gaussian distribution in question. These graphical models provide an underlying structure for the conditional independence relations between variables in the distribution.

The graph learning problem can be summed up as recovering the precision matrix from the sample covariance matrix. Typically, the precision matrix can be learnt by simply inverting the covariance matrix. However, in the high-dimensional setting, a large sample size is required for the inversion to be accurate. To reduce the sample complexity, sparsity was introduced into graphical models so that the topology can be learnt with a smaller sample size, of order $\Omega (k $ log $ p)$. The most common notion of sparsity is a bound on the node degree of the graph. With a degree bound, many regression techniques \citep{MeiBuh:06, RavWaiRasYu:11,RenSunZhaZho:14} like the LASSO can be used to recover the neighborhood of each node, and in doing so, to learn the underlying graph structure. However, while bounding the degree of the graph makes the graph simpler, there are graphs that are simple to learn but are not degree bounded, such as a star graph for example.  Another notion of sparsity was introduced that overcomes the limits of bounded degree graphs and applies to a wider class of graphs
\citep{SohTat:14,AnaTanHuaWil:12}. This sparsity is determined by the measure of separability in a graph, which limits the number of vertex disjoint paths between non-neighboring nodes in the graph. GGMs with separability, quantified by parameter $k$, can also be learnt with  $\Omega (K $ log $ p)$ sample complexity as well.

Besides bounded degree and separability, another type of graph that can be considered to have relatively simple structure is known as graphs with feedback vertex sets (FVSs) \citep{Vaz:04}. A FVS is a node set in a graph whereby the removal of the FVS and any edge in the graph connected to a node in the FVS results in a subgraph of the original which has tree structure. It has been shown that inference on GGMs with FVSs can be done using message passing algorithms \citep{LiuChaAnaWil:12}. In particular, graph learning can also be performed on this type of GGM \citep{LiuWil:13}. In a fully observed GGM, the authors considered the cases where the FVS is known and where the FVS is not known, and provide algorithms for graph learning in both cases. In the event where the FVS is not known, the authors have provided empirical evidence that their algorithm performs reasonably well in identifying the nodes in the FVS.

In our work, we will tackle the problem of identifiability in GGMs. In graph learning algorithms, the sparsity constraint is usually an assumption imposed on the graph, that we treat as prior knowledge of the GGM. If the GGM is known to be sparse, then the GGM can be learnt. However if the GGM is not known to be sparse, the graph learning technique cannot be applied. In this paper, we will do away with the sparsity constraint, and introduce an identifiability algorithm that serves as a sparsity check. Namely, we do not assume a priori that the GGM is sparse. Rather, for any general GGM, we can check using our algorithm to see if the GGM is indeed sparse or not. We do this for degree bounded graphs. In the case of separability however, we introduce a new concept of separability, known as strong separability, that is motivated by the typical notion of separability. We introduce an algorithm that can also determine whether a GGM is strongly separable or not. In the event that the GGM is indeed sparse (bounded degree or strongly separability), we show that our algorithms are able to learn the topology of the GGM as well. We provide theoretical guarantees that the identification and graph learning can be done with $\Omega(k \log p)$ sample complexity, where $k$ is the measure of sparsity and $p$ is the dimension of the GGM.

Also, we extend the notion of a FVS to that of a generalized FVS. A generalized FVS is a node set in the graph such that the removal of the FVS and its connecting edges results in a strongly separable graph. We will show that a strongly separable graph is a natural generalization of a tree. We propose a new algorithm for identifying the nodes that belong to a generalized FVS in a GGM and performs graph recovery on the GGM as well. We provide theoretical guarantees for selecting the FVS in the graph and for learning, and establish that these can be done with $\Omega( (k + \ell) \log p)$ sample complexity, where $k$ captures the measure of strong separability in the FVS removed subgraph, and $\ell$ is the size of the FVS.

The rest of the paper is structured as follows.
In Section \ref{sect:Prelim}, we lay down some preliminary results needed to establish the correctness of our algorithms.
In Section \ref{sect:Degree}, we introduce an algorithm to identify whether a GGM is degree bounded or not.
In Section \ref{sect:StrongKSep}, we discuss the notion of a strongly separable GGM, and propose an algorithm that can determine whether a GGM is strongly separable of not.
In Section \ref{sect:FVS}, we introduce the notion of a generalized FVS, andl extend the above techniques to recover the graph structure of a GGM with generalized FVSs.
We conclude in Section \ref{sect:Conclusion}.

\section{Preliminaries}
\label{sect:Prelim}

Let $\bs{A} \in \mathbb{R}^{p \times p}$ be a matrix, with its index set $\cV = \{1, \ldots p\}$. For $I, J \subset \cV$, $\bs{A}_{IJ}$ is the $|I| \times |J|$ submatrix of $\bs{A}$ corresponding to the formed by the intersection of the rows of $I$ and the columns of $J$, with $\bs{A}_I = \bs{A}_{II}$. If $I^c = \cV \setminus I$, then $\bs{A}$ has the block structure, (with row and column exchanges);
\begin{equation}
\bs{A} = \begin{bmatrix} \bs{A}_{I} & \bs{A}_{I I^c} \\ \bs{A}_{I^cI} & \bs{A}_{I^c} \end{bmatrix}.
\end{equation}
In terms of notation, let $|\bs{A}|$ denote its determinant, $\text{tr}(\bs{A})$ denote its trace, $A^T$ denote its transpose, $\|\bs{A}\|_2$ denote its spectral norm, and  $\lambda_{\max}(\bs{A})$ and $\lambda_{\min}(\bs{A})$ denote its maximum and minimum eigenvalues respectively.  In this paper, we will often refer to the $p$ by $p$ covariance matrix $\bs{\Sigma}$, so we will use the shorthand $\lambda_{\max} = \lambda_{\max} (\bs{\Sigma})$ and $\lambda_{\min} = \lambda_{\min} (\bs{\Sigma})$. 
The Schur complement of $\bs{A}_I$ in $\bs{A}$ is defined by
\begin{equation}
\bs{A}_{I^c \mid I} = \bs{A}_{I^c} - \bs{A}_{I^cI}\bs{A}_I^{-1}\bs{A}_{II^c}.
\end{equation}

In terms of graphs, we consider only graphs without self loops and multiple edges between the same pair of nodes. Let $\cG = (\cV,\cE)$ be an undirected graph, where $\cV = \{1, \ldots, p\}$ is the set of nodes and $\cE$ is the set of edges, where $(u,v) \in \cE$ if and only if an edge exists between nodes $u$ and $v$. Any node connected to a node $u$ by an edge is known as a neighbor of $u$ and the number of distinct neighbors of $u$ is known as the degree of $u$. We denote the set of neighbors of $u$ as $\cN(u)$.
For $I \subseteq \cV$, we denote the induced subgraph on nodes $I$ by $\cG_I$. For two distinct nodes $u$ and $v$, a path of length $t$ from $u$ to $v$ is a series $\{(u, w_1), (w_1, w_2), \ldots, (w_{t-2},w_{t-1}), (w_{t-1},v)\}$ of edges in $\cE$, where $w_1, \ldots, w_{t-1} \in \cV$. The graph $\cG$ is connected if for any distinct nodes $u,v \in \cV$, there is at least one path from $u$ to $v$. Otherwise, the graph $\cG$ is disjoint. A connected component of $\cG$ is a subgraph of $\cG$ that is connected. A disjoint graph can be divided into a number of connected components, where nodes from distinct connected components are not connected by a path.

Let $\bs{X} = (X_1, \ldots, X_p)$ be a multivariate Gaussian distribution with mean $\bs{\mu}$ and covariance matrix $\bs{\Sigma}$. For the rest of this paper, we will only consider zero mean Gaussian distributions, that is, $\bs{\mu} = \bs{0}$. The precision matrix of the Gaussian distribution is $\bs{\Omega} = \bs{\Sigma}^{-1}$. The random variable $\bs{X}$ has the distribution function
\begin{equation}
f_{\bs{X}}(\bs{x}) = \frac{1}{\sqrt{(2\pi)^p|\bs{\Sigma}|}} \exp\left\{ -\frac{1}{2}(\bs{x} - \bs{\mu})^T \bs{\Omega} (\bs{x} - \bs{\mu})  \right\}.
\end{equation}
We denote the independence of $X_u$ and $X_v$ by $X_u \perp X_v$, and the conditional independence of  $X_u$ and $X_v$ given another random variable $X_z$ by $X_u \perp X_v \mid Z$. The precision matrix can be expressed as conditional independence relationships of variables in $\bs{X}$. More precisly, $\bs{\Omega}_{uv} = 0$ if and only if $X_u \perp X_v \mid \bs{X}_{\cV \setminus \{u,v\}}$. 
Given  $S \subseteq \cV$, the conditional distribution of $\bs{X}_{S^c}$ given $\bs{X}_S = \bs{x}_S$ is a multivariate Gaussian distribution with conditional mean $ \bs{\mu}_{S^c} - \bs{\Sigma}_{S^cS} \bs{\Sigma}^{-1}_{S} (\bs{x}_S - \bs{\mu}_S)$ and  conditional covariance matrix 
\begin{equation}
\bs{\Sigma}_{S^c \mid S} =  \bs{\Sigma}_{S^c} - \bs{\Sigma}_{S^cS} \bs{\Sigma}^{-1}_{S} \bs{\Sigma}_{SS^c}.
\label{eq:CondCovMat}
\end{equation} 
Observe that the conditional covariance is the Schur complement of $\bs{\Sigma}_S$ in $\bs{\Sigma}$. It follows then that 
$X_u \perp X_v \mid \bs{X}_S$ if and only if $(\bs{\Sigma}_{S^c \mid S})_{uv} = 0$, so conditional independence relations can be computed from the covariance matrix. The sample complexity in calculating these conditional independence relations depends primary on the size of the set conditioned upon, with the sample complexity being small when the size of the conditioned set is small. 
To denote conditional covariance, we use the notations $\bs{\Sigma}(u,v \mid S)$ and $(\bs{\Sigma}_{S^c \mid S})_{uv}$ interchangeably. 

A Gaussian graphical model therefore is a Gaussian multivariate distribution $\bs{X}$ with a graph  $\cG_{\bs{\Sigma}} = (\cV, \cE)$ associated with it. The node set $\cV$ is the index set of the distribution $\{1, \ldots, p\}$, and $(u,v)$ is in $\cE$ if and only if $\bs{\Omega}_{uv} \neq 0$. The graph   $\cG_{\bs{\Sigma}} $ is known as a precision or concentration graph. For simplicity, we will mostly refer to the precision graph as $\cG$ unless there is ambiguity.

\subsection{Sample Analysis}
\label{sect:SampleCovariance}

Let $\bs{x}^{(1)}, \ldots, \bs{x}^{(n)} \in \mathbb{R}^p$ be $n$ i.i.d. samples of the random variable $\bs{X}$ with distribution $\cN(\bs{0},\bs{\Sigma})$.  The scatter matrix $\bs{S}$ is defined as
\begin{equation}
\bs{S} = \sum_{i=1}^n \bs{x}^{(i)} (\bs{x}^{(i)})^T.
\end{equation}
The sample covariance matrix determined by these $n$ samples is defined as
\begin{equation}
\widehat{\bs{\Sigma}} = \frac{1}{n} \bs{S}.
\end{equation}
In determining the sample conditional covariances, we will make use of the scatter matrix $\bs{S}$ instead of $\widehat{\bs{\Sigma}} $.
Let $u$ and $v$ be distinct elements of $\cV$ and let $S \subseteq \cV \setminus \{u,v\}$. The sample conditional covariance of $X_u$ and $X_v$ given $\bs{X}_S$ is denoted by
\begin{equation}
\widehat{\bs{\Sigma}}(u,v \mid S) = \frac{1}{n - |S|} \left(\bs{S}_{uv} - \bs{S}_{uS}\bs{S}_S^{-1}\bs{S}_{Sv} \right).
\label{eq:SampleConditional}
\end{equation}
In our algorithms, we usually have to decide whether a conditional independence relation holds. We have to determine whether $X_u \perp X_v \mid \bs{X}_S$ or $X_u \not\perp X_v \mid \bs{X}_S$. To do so with the sample covariance matrix, we need to define a conditional independence threshold $\alpha >0$, such that if 
\begin{equation}
|\widehat{\bs{\Sigma}}(u,v \mid S)| < \alpha,
\end{equation} 
we will decide that $X_u \perp X_v \mid \bs{X}_S$. Otherwise, we decide that $X_u \not\perp X_v \mid \bs{X}_S$. In our analysis, $\alpha$ will scale depending on $p, n$ and $|S|$.

\subsection{Faithful Conditional Independence Relationships}

In this paper, we will use conditional covariances and independences to determine graph relationships. As mentioned previously, conditional independence is closely related to separability in graphs. The example we used before was the local Markov property. There are two other Markov properties that hold in Markov random fields. The first is the pairwise Markov property that states $X_u \perp X_v \mid \bs{X}_{\cV \setminus \{u,v\}}$. The second is the global Markov property, which states that if $S$ is a vertex separator of $u$ and $v$, then $X_u \perp X_v \mid \bs{X}_S$. 

The Markov properties describe how graph structures imply conditional independence relationships. However, in graph learning, we want to deduce graph structure from conditional independence relations. To do so, we will borrow some results regarding the faithfulness of conditional independence relations \citep{SohTat:14}, which can also be described as when the converse of the global Markov property holds.

\begin{definition}
Let $\bs{X}$ be a Gaussian graphical model. A conditional independence relation $X_u \perp X_v \mid \bs{X}_S$ is said to be faithful if $S$ is a vertex separator of $u$ and $v$ in the precision graph $\cG$. Otherwise $X_u \perp X_v \mid \bs{X}_S$ is unfaithul.
\end{definition}

Using other conditional relationships in a graph, we can check whether a particular conditional independence relationship is faithful or not. We will use this result later in our graph learning algorithm. Given a conditioning node set $S \subset \cV$, let the graph $\bar{\cG}^{S^c} =(S^c, \bar{\cE}^{S^c})$ be defined by $S^c = \cV \setminus S$, and $(i,j) \in \bar{\cE}$ if and only if $X_i \not\perp X_j \mid \bs{X}_S$. By observing the connected components of $\bar{\cG}$, we have the following proposition. 

\begin{proposition}
\label{Prop:Faithfulness}
Let $\bs{X}$ be a Gaussian graphical model with precision graph $\cG$. A conditional independence relation $X_u \perp X_v \mid \bs{X}_S$ is faithful if and only if $u$ and $v$ are in separate connected components of $\bar{\cG}^{S^c}$. Also any two nodes that are in separate connected components  in $\bar{\cG}^{S^c}$ are not connected by an edge in $\cG$.
\end{proposition}


\section{Degree Bounded GGMs}
\label{sect:Degree}

In this section, we will look at a specific class of GGMs known as degree bounded GGMs, and we will use degree bounded GGMs to illustrate the identifiability problem and technique. Many high-dimensional graph learning techniques typically assume that the GGMs have a $k$-degree bounded graph structure, making it possible to learn the topology of the GMM often with a sample complexity of $\Omega( k \log p)$, where $p$ is the dimension of the GGM. However, now suppose we do not know a priori that the GGM has a degree bounded graph. In order for these graph learning techniques to work, we first need to identify whether the GGM has a degree bounded graph or not. Therefore, we will introduce an algorithm that can do so for any general GGM. In the process, if the underlying graph is found to be degree bounded, we will end up learning the graph structure as a result.

We begin by defining a $k$-degree bounded graph.
\begin{definition}
A graph $\cG = (\cV,\cE)$ is said to be a $k$-degree bounded graph if $|\cN(u)| \leq k$ for all $u \in \cV$.
\end{definition}
In order to describe the identification procedure for a GGM, we have to define certain vertex sets.
For a node $u \in \cV$, we define the set
\begin{align}
S^{\text{deg}}_k(u)  =  & \: \{ S \subset \cV \setminus \{u\} : |S| \leq k, \nonumber \\   & \:
 X_u \perp \bs{X}_{\cV\setminus S \cup \{u\}} \mid \bs{X}_S \}.
\end{align}
In our algorithm, for every node $u$, we will check to see if the neighbor set of $u$, $\cN(u)$, is not greater than $k$. If we condition on the variable $\bs{X}_{\cN(u)}$, then $X_u$ will be conditionally independent of the variables $\bs{X}_{\cV \setminus (\{u\} \cup \cN(u))}$. The algorithm searches for a set $S$ of size $k$ that contains the $\cN(u)$, which is possible if $\cN(u)$ has size at most $k$. If such a set can be found for every node $u \in \cV$, then Algorithm \ref{Algo:Deg} will decide that the GGM is $k$-degree bounded. In this way, Algorithm \ref{Algo:Deg} differentiates between GGMs that are $k$-degree bounded and those that are not. Also, if any set $S$ that separates $u$ from the rest of the nodes $\cV \setminus (S \cup \{u\})$, then $u$ is not connected by an edge to any node in $\cV \setminus (S \cup \{u\})$. This means that the neighbors of $u$ must be in $S$, and the taking the intersection of all possible $S$ will give us the neighbor set of $u$. It follows then that if the GGM is $k$-degree bounded, we can also deduce exactly the topology of $\cG$, which gives us the following theorem.

\begin{algorithm}[tb]
   \caption{Identifying a $k$-degree bounded GGM $\bs{X}$ }
   \label{Algo:Deg}
\begin{algorithmic}
   \STATE {\bfseries Input:}Covariance matrix $\bs{\Sigma}$ and degree bound $k$.
   \IF{ For all $u \in \cV$, $S^{\text{deg}}_k(u) \neq \phi$}
   \STATE Output $\cG$ as $k$-degree bounded.
   \STATE Output $\cN (u) = \bigcap_{S \in S^{\text{deg}}_k(u)} S$. 
   \ELSE
   \STATE Output  $\cG$ as not $k$-degree bounded.
   \ENDIF
\end{algorithmic}
\end{algorithm}

\begin{theorem}
Given a GGM $\bs{X}$  with covariance matrix $\bs{\Sigma}$, Algorithm \ref{Algo:Deg} correctly identifies whether the  underlying graph $\cG$ of $\bs{X}$ is a $k$-degree bounded graph or not. In the case where it is a $k$-degree bounded graph, Algorithm \ref{Algo:Deg} also correct ouputs the neighbor set of all nodes $u \in \cV$, and consequently, recovers the structure of the graph $\cG$ of the GGM.
\label{Theo:DegreeIdent}
\end{theorem}

\textit{Proof.}
Suppose a GMM is $k$-degree bounded. For any node $u \in \cV$, its set of neighbors does not have size more than $k$, and so any size $k$ node set that contains the neighbor set is separates $u$ from the rest of the nodes. By Proposition \ref{Prop:Faithfulness}, this means the $S^{\text{deg}}_k(u)$ contains at least that particular size $k$ node set that separates $u$ from the rest of the nodes, and is thus non-empty. Therefore Algorithm \ref{Algo:Deg} correctly outputs the GGM as a $k$-degree bounded graph. Now suppose that the GMM is not $k$-degree bounded. Then there exists a node $u$ whereby its neighbor set is larger than $k$. By Proposition  \ref{Prop:Faithfulness} , any sets that separates $u$ from the rest of the nodes must contain its neighbor set, therefore $S^{\text{deg}}_k(u)$ is empty, therefore Algorithm \ref{Algo:Deg} will correct output the GGM as being not a $k$-degree bounded graph as well.
\QEDA

\begin{corollary}
Let $\bs{X}$ be a GMM with covariance matrix $\bs{\Sigma}$. Let $\bs{S}$ be used according to \eqref{eq:SampleConditional} to determine the sample conditional covariances in the procedure outlined in Algorithm \ref{Algo:Deg}, instead of the true covariance matrix $\bs{\Sigma}$. 
Let $\beta \leq \min_{|S| = k , X_u \not\perp X_v \mid \bs{X}_S} |\bs{\Sigma} (u,v \mid S)|$, and let $\alpha = \beta /2$. Let $M_0$ be the event that Algorithm \ref{Algo:Deg} correctly identifies whether $\cG$ is a $k$-degree bounded graph or not.
Then,
\begin{equation}
\P\left( M_0 \right) \geq 1 - \epsilon ,
\end{equation}
and
\begin{equation}
\P\left(\hat{\cE} = \cE \mid \cG \:
 \text{is $k$-degree bounded} \right) \geq 1 - \epsilon ,
\end{equation}
for $n = \Omega \left( \frac{\lambda_{\max}^2 + \beta\lambda_{\max}}{\beta^{2}}( k \log p + \log (\epsilon^{-1})) \right)$.
\label{Corr:DegreeSample}
\end{corollary}

The proof of this theorem follows from the theoretical results in the work on separable graphs \citep{AnaTanHuaWil:12} .
The main idea is that with correcting scaling of the sample complexity, we will have high probability of accurately identifying the conditional independence relations in the graph. The scaling therefore, depends on the size of the set $S$ that is conditioned on in these conditional independence relationships, which is upper bounded by $k$, hence giving us the required sample complexity result.

This algorithm considers the node sets in $\cV$ that have at most size $k$ in order to capture the neighbor set, so it has a computational runtime of $O(p^{k+1})$. However, when the neighbor set for one node $u$ is identified, immediately all the other non-edges between $u$ and its non-neighbors are known. Now, for any of the neighbors of $u$, the runtime to find its neighbors is reduced by a factor of $p$ since it is already known that $u$ is one of its neighbors. In this way, the actual runtime can be reduced by storing the edges and non-edges of the graph as we learn them.


\section{Strongly $k$-separable GGMs}
\label{sect:StrongKSep}

In this section, we will present an identifiability algorithm for separable graphs. Separability is an emerging concept in the area of statistical graph learning. Naturally so, since many problem in graph learning involves recovering edges or entries in the precision matrix, which are parameters quantifying the interaction between pairs of nodes. Separability captures the number of vertex disjoint paths betwen pairs of nodes, and as a result, shows itself to be a better measure of the relationship between pairs of nodes compared to more localized properties of the graph such the node degree. Also, any graph with bounded degree is a separable graph, and as such, graph learnings for separable graphs can be applied to bounded degree graphs as well. In other words, a broader class of graphs are learned by looking at separable graphs. 

\begin{definition}
Let $u$ and $v$ be two non-neighboring nodes in a graph $\cG$. Suppose there exists a node set $S \in \cV \setminus \{u,v\}$ with $|S| = k$ such that every path from $u$ to $v$ has to pass through some node in $S$. Then $u$ and $v$ are said to be $k$-separable. A graph $\cG$ is said to be $k$-separable if every pair of non-neighboring nodes in the graph is $k$-separable.
\end{definition}

A $k$-separable GGM can be learnt if its $k$-separability is known or assumed from the start. However, $k$-separable GGMs could potentially be hard to identify due to the little restriction $k$-separability places on the density of a graph, since we want to do so with low sample complexity as well. In $k$-separable graphs, cliques of arbitrary size could be present. A complete graph, for exmaple, is a $k$-separable graph for any $k$. Therefore, any identifiability algorithm for $k$-separable GGMs must be able to distinguish between dense graphs, such as between the complete graph and the complete graph with a single edge missing. 

The sparsity comes into play in the identifiability problem mainly because of the low sample complexity we are trying to achieve. In the spirit of making the graph less dense, we will extend the idea of $k$-separability, and generalize the definition to pairs of nodes that are neighbors as well. It turns out that in doing so, this new notion of separability is actually identifiable. 

\begin{definition}
Let $u$ and $v$ be two nodes in a graph $\cG$. Let $\delta_{uv}$ be equal to $0$ if $u$ and $v$ are not connected by an edge, and let it be equal to $1$ otherwise. Suppose there exists a node set $S \in \cV \setminus \{u,v\}$ with $|S| = k -\delta_{uv}$ ,such that every path from $u$ to $v$, excluding the edge $(u,v)$ if it exists, has to pass through some node in $S$. Then $u$ and $v$ are said to be strongly $k$-separable. A graph $\cG$ is said to be strongly $k$-separable if every pair of nodes in the graph is strongly $k$-separable.
\end{definition}

By the above defintion, two pair of non-neighboring nodes are strongly $k$-separable if and only if they are $k$-separable. Thus, strong separability can be seen as an extension of the typical notion of separability extended to neighboring nodes as well. Here are some examples of strongly $k$-separable graphs.

\begin{lemma}
A connected strongly $1$-separable graph is a tree. A connected strongly $2$-separable graph is a series of rings and trees iteratively connected nodewise. 
\label{Lemm:1Sep}
\end{lemma}

\textit{Proof.}
By definition, no cycles exist in a strongly $1$-separable graph, since any pair of nodes that form an edge in the cycle is connected by another path other than their edge and thus is not strongly $1$-separable. Therefore, a connected strongly $1$-separable graph must be a tree. 

Any strongly $2$-separable graph is a $2$-separable graph. Thus, it can be expressed as a series of rings, trees or cliques iteratively connected by merging a common node \citep{CicMil:12}. Since a strongly $2$-separable clique is an edge or a triangle, which is also a cycle, we therefore have our result.
\QEDA

\begin{lemma}
A $k$-degree bounded graph is a strongly $k$-separable graph.
\end{lemma}

\textit{Proof.}
Let $\cG$ be a $k$-degree bounded graph. For any pair of non-neighboring nodes $u$ and $v$, the neighbor set $\cN(u)$ is a  vertex separator set of $u$ and $v$ of size $k$, since the graph is degree bounded. For any node $u$ and a neighboring pair $v$, the neighbor set $\cN(u) \setminus \{v\}$ separates $u$ has size $k-1$, and all paths from $u$ to $v$ must pass through $\cN(u) \setminus \{v\}$. This means that $u$ and $v$ are strongly $k$-separable whether they are connected by an edge or not, and therefore $\cG$ is strongly $k$-separable.  
\QEDA

\begin{example}
We construct a graph $\cG$ as follows, where $\cV = \{1, \ldots, p\}$. Connect $1$ and $2$, and let nodes the neighbor set of $N(u)$ be $\{1,2\}$ for $u \geq 3$. Then for any $k < p-2$, we have that $\cG$ is $k$-separable but not strongly $k$-separable. 
\end{example}

With this generalization of separability, these strongly $k$-separable GGMs are now identifiable. Given a GGM without prior knowledge of its structure, we can determine whether the GGM is strongly $k$-separable or not, with $\Omega (k \log p)$ sample complexity. Tha main idea is to check whether each node pair is strongly $k$-separable. There are two cases involved, the first is where the nodes in the pair are not neighbors, and the second is where they are. If the node pair are non-neighbors and they are strongly $k$-separable, then this can be identified via the work in learning $k$-separable graphs \citep{SohTat:14}. Thus, the only other case where the node pair can be strongly $k$-separable is if they were neighbors. We will now introduce a method to test if this is true and assimilate this into our identification algorithm.

To determine the case where $(u,v)$ is not an edge and is strongly $k$-separable, we will define the set
\begin{align}
S^{\text{sep}}_k(u,v) = & \: \{S \in \cV \setminus \{u,v\} : |S| \leq k, \nonumber \\
& \: \: \Sigma(u,v \mid S) = 0 \: \text{and is faithful}\}. 
\end{align}
This set captures all the possible node sets of size $k-1$ that separate $u$ and $v$.

Suppose a node pair $(u,v)$ is strongly $k$-separable and is an edge in the graph. Let the subgraph $\cG_{-(u,v)}$ to be the graph $\cG$ with edge $(u,v)$ removed. 
There is a node set $S \subset \cV \setminus \{u,v\}$, $|S| \leq k-1$, such that the removal of edge $(u,v)$ from $\cG$ results in a graph where $S$ separates nodes $u$ and $v$.  This separation property is exactly where we can make use of Proposition \ref{Prop:Faithfulness} to identify if a node pair is strongly $k$-separable.
If we can remove the edge $(u,v)$ from the graph, we can then use conditional independence relations to find a node separator $S$ of $u$ and $v$ in the resultant graph $\cG_{-(u,v)}$. Of course, we cannot simply remove this edge simply based on the covariance matrix. 
We also cannot simply condition on $\bs{X}_S$, because, if $(u,v)$ is in $\cE$, we will have the condition dependence relation $X_u \not\perp X_v \mid \bs{X}_S$. We could try remove the influence of the edge $(u,v)$ in the graph by conditioning on $\bs{X}_{S \cup \{u,v\}}$, however this does not ensure that $u$ and $v$ are separated by $S$ in $\cG_{-(u,v)}$.

To overcome this problem, we condition on both $\bs{X}_{S \cup \{u\}}$ and $\bs{X}_{S \cup \{v\}}$. We use the conditional independence relations given these random variables to deduce that $S$ is a node separator of $u$ and $v$ in $\cG_{-(u,v)}$. Running through node subsets $S \subset \cV \setminus \{u,v\}$ of size $k-1$, we first condition on $\bs{X}_{S \cup \{v\}}$ to see how $S$ separates $\cG_{S^c \setminus \{v\}}$. We then condition on $\bs{X}_{S \cup \{u\}}$ to see how $S$ separates $\cG_{S^c \setminus \{u\}}$. Using these two pieces of information, and paying attention to the connected components that arise in both cases, we can infer whether $S$ separates $u$ and $v$ in $\cG_{-(u,v)}$.

For any subset $S \subset \cV$, we define the graph $\bar{\cG}^{S^c} = (S^c, \bar{\cE}^{S^c})$, where $(u,v) \in \bar{\cE}^{S^c}$ if and only if $\bs{X}_u \not\perp X_v \mid \bs{X}_S$. For a node $ h \in S^c$, let the connected node set component of $\bar{\cG}^{S^c}$ containing $h$ be denoted by $\bar{U}_{S^c} (h)$.

For any node $u \in \cV$, we denote the set
\begin{equation}
\Gamma_{(u,v)}^k = \{ S \subset \cV  \setminus \{u,v\}: |S| \leq k-1\}.
\end{equation}
of all possible node subsets $S$ of size $k-1$ in $\cV \setminus \{u,v\}$. 
We define a subset of this set, which is
\begin{align}
\Gamma^k_{u \mid v} = & \: \{ S \in \Gamma^k_{(u,v)}: \: \exists \: h \in S^c \setminus  \{u,v\}\:\: \text{s.t.}\nonumber \\ 
& \:\: \bs{\Sigma} (u,h \mid S+v) = 0 \:\: \text{and is faithful} \},
\end{align}
where $S+v = S \cup \{v\}$, and the faithfulness of the relation $\bs{\Sigma} (u,h \mid S+v) = 0$ is determined by Proposition \ref{Prop:Faithfulness}. This quantity encompasses the different sets $S$ such that $\cG_{S \cup \{v\}}$ is a disjoint graph. However, this set does not subsume all possible $S$ that separate $u$ and $v$ in $\cG_{-(u,v)}$. To include all such possible node sets $S$, we specify a subset of $\Gamma_{u \mid v}$, namely,
\begin{align}
\Psi^k_{u \mid v} =& \: \{ S \in \Gamma^k_{(u,v)}: \:  \bs{\Sigma} (u,h \mid S +j) = 0,\nonumber \\
&  \: \: \forall h \in S^c \setminus  \{u,v\}\}.
\end{align}
These sets cater specifically to the case where $S$ neighbor separates $u$ and $v$, but $u$ has only one neighbor, $v$, in $\cG_{S^c}$. Finally,  let 
\begin{align}
\Lambda^k_1(u,v)  = & \: \{S \in \Gamma_{u \mid v} \cap \Gamma_{v \mid u} : \nonumber \\
& \: \:  \bar{U}_{S^c \setminus \{u\}} (v)  \subseteq \left( S^c \setminus \bar{U}_{S^c \setminus \{v\}} (u) \right) \}.
\end{align}
Also, define
\begin{equation}
\Lambda^k_2(u,v) = \Psi_{v \mid u},
\end{equation}
and
\begin{equation}
\Lambda^k_3(u,v) = \Psi_{u \mid v}.
\end{equation}
Finally, let
\begin{equation}
\Lambda^k_0(u,v) = \Lambda^k_1(u,v) \cup \Lambda^k_2(u,v) \cup \Lambda^k_3(u,v).
\end{equation}
If $(u,v)$ is a pair of neighbor nodes that is strongly $k$-separable, then the set $\Lambda_0(u,v)$ is non-empty. In this way we can identify whether the node pair is a neighbor pair that is strongly $k$-separable. 

Our algorithm therefore aims to check for each node pair if they are indeed strongly $k$-separable. If this doesn't hold for any node pair, then the algorithm will infer that the graph is not strongly $k$-separable. In the case where this holds for every node pair, then the algorithm will output that the graph is strongly $k$-separable. In the process, we get additional information in this case, since the algorithm will also tell us if $u$ and $v$ are connected by edge or not. Thus, when the graph is strongly $k$-separable, we can learn the graph topology as well. This leads to the following theorem.

\begin{algorithm}[tb]
   \caption{Identifying a strongly $k$-separable GGM $\bs{X}$ }
   \label{Algo:Strong}
\begin{algorithmic}
   \STATE {\bfseries Input:}Covariance matrix $\bs{\Sigma}$ and parameter $k$.
   \STATE {\bfseries Initialize:} Pair Set $P = \phi$.	
   \FOR{$u,v \in \cV$, $u > v$}
   \IF{$S^{\text{sep}}_k(u,v) \neq \phi$}
   \STATE Add $(u,v)$ to set $P$.
   \STATE Output $(u,v)$ as non-neighbors.
   \ELSIF{$\Lambda_0^k(u,v) \neq \phi$}
     	\STATE Add $(u,v)$ to set $P$.
	\STATE Output $(u,v)$ as neighbors.  
	\ENDIF
   \ENDFOR
   \IF{$P = \{(u,v) \in \cV \times \cV : u > v\}$}
   \STATE Output $\cG$ as being strongly $k$-separable.
   \STATE Output $\hat{\cE}$.
   \ELSE
   \STATE Output $\cG$ as being not strongly $k$-separable.
   \ENDIF
\end{algorithmic}
\end{algorithm}

\begin{theorem}
Given a GGM $\bs{X}$  with covariance matrix $\bs{\Sigma}$, Algorithm \ref{Algo:Strong} correctly identifies whether the  underlying graph $\cG$ of $\bs{X}$ is a strongly $k$-separable graph or not. In the case where it is a strongly $k$-separable graph, Algorithm \ref{Algo:Strong} also correctly recovers the structure of the graph $\cG$ of the GGM.
\label{Theo:StrongIdent}
\end{theorem}

\textit{Proof.}
A pair of nodes $(u,v)$ is non-neighboring and is $k$-separable if and only if  $S^{\text{sep}}_k(u,v)$ is non-empty \citep{SohTat:14}. We will first show a similar relation for neighboring nodes, namely, if $u$ and $v$ are neighbors, then they are strongly $k$-separable if and only if $\Lambda_0^k(u,v)$ is non-empty.

Suppose a pair of neighbor nodes $\{u,v\}$ is strongly $K$-separable. Then $\Lambda_0^k(u,v)$ must be non-empty, since there is a set $S$ of at most size $k-1$ such that any path from $u$ to $v$ must contain a node in $S$, excluding the edge between $u$ and $v$. Thus, $S \in \Lambda_0^k(u,v)$. Now suppose $\Lambda_0^k(u,v)$ is non-empty, then any node set $S \in \Lambda_0^k(u,v)$ has the follow properties: The graph $\cG_{S^c}$ can be partitioned into four nodes sets $\{u\}$, $\{v\}$, $R_u$ and $R_v$, with $R_u =  \bar{U}_{S^c \setminus \{v\}} (u)$ and $R_v = S^c \setminus (\{u,v\} \cup  \bar{U}_{S^c \setminus \{v\}} (u))$. These four sets have the property that there are no edges connecting any node from $R_u$ to any node in $R_v \cup \{v\}$ and there are no edges connecting any nodes from $R_v$ to $R_u \cup \{u\}$, by Proposition \ref{Prop:Faithfulness}. This means that the only possible edge between $R_u \cup \{u\}$ and $R_v \cup \{v\}$ is the edge between $u$ and $v$. Therefore, $S$ separates $u$ and $v$ in $\cG_{-(u,v)}$, so $u$ and $v$ are strongly $k$-separable.

Let the GGM $\bs{X}$ be strongly $k$-separable. Then every node pair in $\cG$ is either a pair of non-neighbors or a pair of neighbors. If the pair is not connected by an edge, then $S^{\text{sep}}_k(u,v)$ is non-empty. If the pair is connected by an edge,  then$\Lambda_0^k(u,v)$ is non-empty. Therefore, Algorithm \ref{Algo:Strong} will output the graph correctly as a strongly $k$-separable graph. Next, suppose the GGM $\bs{X}$ is not strongly $k$-separable. Then there exists a pair of nodes that is not strongly $k$-separable. There are two cases. If $u$ and $v$ are non-neighbors, then clearly $S^{\text{sep}}_k(u,v)$ must be empty. Also, if  $\Lambda_0^k(u,v)$ is non-empty, then $u$ and $v$ would be $(k-1)$-separable, which is a contradiction. Thus $\Lambda_0^k(u,v)$  must be empty. In the second case, $u$ and $v$ are neighbors. Then $\Lambda_0^k(u,v)$  must be empty. Also, since $u$ and $v$ are connected by an edge,  $S^{\text{sep}}_k(u,v)$ must be empty as well as $u$ and $v$ cannot be separated by an node set. Therefore, $(u,v)$ will not be placed in $P$ by Algorithm \ref{Algo:Strong} and so, Algorithm \ref{Algo:Strong} will output $\cG$ as not being strongly $k$-separable. Consequently, Algorithm \ref{Algo:Strong} correctly identifies whether a GGM is strongly $k$-separable or not.
\QEDA

\begin{corollary}
Let $\bs{X}$ be a GMM with covariance matrix $\bs{\Sigma}$. Let $\bs{S}$ be used according to \eqref{eq:SampleConditional} to determine the sample conditional covariances in the procedure outlined in Algorithm \ref{Algo:Strong}, instead of the true covariance matrix $\bs{\Sigma}$. 
Let $\beta \leq \min_{|S| = k, X_u \not\perp X_v \mid \bs{X}_S} |\bs{\Sigma} (u,v \mid S)|$, and let $\alpha = \beta /2$. Let $M_1$ be the event that Algorithm \ref{Algo:Strong} correctly identifies whether $\cG$ is a strongly $k$-separable graph or not.
Then,
\begin{equation}
\P\left( M_1 \right) \geq 1 - \epsilon ,
\end{equation}
and
\begin{equation}
\P\left(\hat{\cE} = \cE \mid \cG \:
 \text{is strongly $k$-separable} \right) \geq 1 - \epsilon ,
\end{equation}
for $n = \Omega \left( \frac{\lambda_{\max}^2 + \beta\lambda_{\max}}{\beta^{2}}( k \log p + \log (\epsilon^{-1})) \right)$.
\label{Corr:StrongSample}
\end{corollary}
Just as in the degree bounded case, the sample complexity results follows directly from the work on separable graphs \citep{AnaTanHuaWil:12}, so we omit the proof.

Algorithm \ref{Algo:Strong} looks through the possible separator sets for each node pair $(u,v)$ to determine the strong separability of the pairs. The computational complexity for this algorithm is $O(p^{k+4})$. The runtime can be further reduced because many redundant steps are included in the algorithm. Whenever the algorithm checks for a separator set for $u$ and $v$, by Proposition \ref{Prop:Faithfulness}, it will also output many non-edges in the graph. Thus by checking for one pair of nodes, we can learn more about other parts of the graph as well.


\section{Learning GGMs with Generalized FVSs}
\label{sect:FVS}

In this section, we will introduce a novel algorithm for learning GGMs with generalized feedback vertex sets (FVSs). Graphs with simple structure are usually easier to learn, in the sense that topological recovery or inference can be done with less sample or computational complexity. It is also for this reason that sparsity constraints like degree bounds or separability are assumed. However, there are graphs that do not adhere to these sparsity guidelines but are mostly simple, where the graph becomes simple if we remove a small number of nodes and the edges connected to them. Graphs with feedback vertex sets are a good example of this.

\begin{definition} 
Let $\cG$ be a graph. A feedback vertex set (FVS) of $\cG$ is a node set $F \subset \cV$ such that the induced subgraph $\cG_{\cV \setminus F}$ is a tree graph.
\end{definition}

This is the typical setting of a FVS. When the nodes in the FVS are removed, along with the edges that connect to them, the resultant graph is a tree graph. In a GGM setting, this means that the conditional distribution $\bs{X}_{\cV \setminus F} \mid \bs{X}_F$ has a precision graph that is a tree structure. Therefore, algorithms designed to learn tree GGMs can be applied to the conditional distribution $\bs{X}_{\cV \setminus F} \mid \bs{X}_F$. In the case where we know exactly which nodes that belong to the FVS, we can proceed to learn the rest of the graph through conditioning on the variables corresponding to the FVS. However, the more challenging problem is that of learning the graph while not knowing where the FVS is in the graph. To do so we have to identify which nodes are in the FVS. 

We will introduce a technique to learn the location of the FVSs in a GGM, and we will do so for a more general class of graphs. By Lemma \ref{Lemm:1Sep}, a strongly $1$-separable graph is a tree. Therefore we can think of a FVS as the removal of a set of nodes that result in a strongly $1$-separable graph. More generally, we can define a generalized FVS whereby the removal of the generalized FVS results in a strongly $k$-separable graph. We provide a formal defintion.

\begin{definition} 
Let $\cG$ be a graph. A $k$-generalized FVS of  size $\ell$ of the graph $\cG$ is a node set $F \subset \cV$, $|F| \leq \ell$ such that the induced subgraph $\cG_{\cV \setminus F}$ is a strongly $k$-separable graph.
\end{definition}

When $k = 1$, the $k$-generalized FVS reduces to the typical FVS. Therefore, the $k$-generalized FVS serves as a generalization of the FVS. Many graphs with generalized FVS do not satisfy the typical sparsity constraints, as such previous graph learning techniques cannot be applied to learn GGMs that contain generalized FVS. In fact, just as separable graphs can be treated as  a generalization of degree bounded graphs, so graphs with generalized FVS can be seen as a generalization of strongly $k$-separable graphs. 

\begin{lemma}
Any strongly $k$-separable graph contains a $k$-generalized FVS of arbitrary size.
\end{lemma}

\textit{Proof.} The removal of any number of edges of the graph does not increase the connectivity or the number of disjoint paths between nodes. Therefore, any node set in the graph is an FVS, the removal of the node set and edges connected to it will preserve the strong separability of the graph.
\QEDA

However, a graph with a generalized FVS may not be strongly $k$-separable. The construction in the following example demonstrates this. 

\begin{example}
Let $\cG' = (\cV', \cE')$ be a tree. We construct $\cG$ by adding to two nodes $u$ and $v$ to $\cG'$, and we connect $u$ and $v$ to every node $\cG'$ by an edge. Let $u$ and $v$ have no edge connecting them. Then in $\cG$, the two nodes $u$ and $v$ are not strongly $k$-separable for any $k < |cV'|$ since the number of vertex disjoint paths between $u$ and $v$ is equal to $|\cV'|$. Then $\cG$ has a $1$-generalized FVS of size $2$. In the same way, an generalized FVS set can be added to a strongly $k$-separable graph so that the resultant graph is not strongly $k$-separable.
\end{example}

Therefore, learning a GGM with a generalized FVS is more general than learning a strongly separable graph. In learning degree bounded graphs or separable graphs, with the sparsity assumption we typically do not need to use identifiability algorithms introduced in previous sections to learn the graph. However, in the case of graphs with generalized FVSs, our learning algorithm still requires us to make use of identifiability techniques. For example, even though we assume the underlying graph of the GGM contains a FVS $F$, we still need to identify whether the induced subgraph  $\cG_{\cV \setminus F}$ is a tree, in order to properly determine which nodes are in $F$. This highlights the importance of the ability to identify certain graph properties without prior assumptions.

We will now introduce an algorithm to identify whether a graph contains a $k$-generalized FVS of size $\ell$, and in the process, discover which nodes belong to the $k$-generalized vertex set. To describe the algorithm, we require the following set to be defined, namely,
\begin{equation} 
S_F = \{S \in \cV : F \subseteq S\}.
\end{equation}
We will also show that this algorithm can learn the structure of the induced subgraph $\cG_{\cV \setminus F}$, where $F$ is a $k$-generalized FVS of size $\ell$.

\begin{algorithm}[tb]
   \caption{Identifying a GGM $\bs{X}$ with a $k$-generalized FVS of size $\ell$}
   \label{Algo:FVS}
\begin{algorithmic}
   \STATE {\bfseries Input:}Covariance matrix $\bs{\Sigma}$ and parameters $k,\ell$.
\FOR{$F \subset \cV, |F| = \ell$}
   \STATE {\bfseries Initialize:} Pair Set $P = \phi$.	
   \FOR{$u,v \in \cV$, $u > v$}
   \IF{$S^{\text{sep}}_k(u,v) \cap S_F  \neq \phi $} 
   \STATE Add $(u,v)$ to set $P$.
   \STATE Output $(u,v)$ as non-neighbors.
   \ELSIF{$\Lambda_0^k(u,v) \cap S_F \neq \phi$}
     	\STATE Add $(u,v)$ to set $P$.
	\STATE Output $(u,v)$ as neighbors.  
	\ENDIF
   \ENDFOR
   \IF{$P = \{(u,v) \in \cV  \times \cV  : u,v \notin F, u > v\}$}
   \STATE Output $F$ as a $k$-generalized FVS.
   \STATE Output $\hat{\cE}$ of $\cG_{\cV \setminus F}$.
   \ELSE
   \STATE Output $F$ as not a $k$-generalized FVS.
   \ENDIF
\ENDFOR
\end{algorithmic}
\end{algorithm}

\begin{theorem}
Given a GGM $\bs{X}$  with covariance matrix $\bs{\Sigma}$, Algorithm \ref{Algo:FVS} correctly identifies all the $k$-generalized FVS  of size $\ell$ in the graph $\cG$. In the case where a $k$-generalized FVS of size $\ell$ exists in the graph, Algorithm \ref{Algo:FVS} also correctly recovers the structure of the induced subgraph $\cG_{\cV\setminus F}$ of the GGM.
\label{Theo:FVSIdent}
\end{theorem}

\textit{Proof.}
Suppose $F$ is a $k$-generalized FVS of size $\ell$. Then $\cG_{\cV \setminus F}$ is a strongly $k$-separable graph. By Theorem \ref{Theo:StrongIdent}, Algorithm \ref{Algo:FVS} will identify the conditional Gaussian distribution $\bs{X}_{\cV \setminus F} \mid \bs{X}_F$ as having a strongly $k$-separable underlying graph. Also, by Theorem \ref{Theo:StrongIdent}, the identification of  $\bs{X}_{\cV \setminus F} \mid \bs{X}_F$  as a strongly $k$-separable  graph will output its underlying topology as well. Therefore, Algorithm \ref{Algo:FVS} will output $F$ correctly as a $k$-generalized FVS of size $\ell$. If $F$ is not a $k$-generalized FVS of size $\ell$, then by Theorem \ref{Theo:StrongIdent}, Algorithm \ref{Algo:FVS} will identify  $\bs{X}_{\cV \setminus F} \mid \bs{X}_F$  as not having a strongly $k$-separable graph structure, and so will identify F correctly as not being a $k$-generalized FVS of size $\ell$. In this way, Algorithm \ref{Algo:FVS} will correctly identify all the FVSs in the GGM. This concludes the proof.
\QEDA

To describe the sample complexity result, for any node set $F\subseteq \cV$, we define $F^c = \cV \setminus F$. We define the subgraph $\cG_{F^c} = (F^c \cE_{F^c})$, so $\cE_{F^c}$ describes the edges that are between the nodes in $F^c$. 

\begin{corollary}
Let $\bs{X}$ be a GMM with covariance matrix $\bs{\Sigma}$. Let $\bs{S}$ be used according to \eqref{eq:SampleConditional} to determine the sample conditional covariances in the procedure outlined in Algorithm \ref{Algo:FVS}, instead of the true covariance matrix $\bs{\Sigma}$. 
Let $\beta \leq \min_{|S| = k + \ell, X_u \not\perp X_v \mid \bs{X}_S} |\bs{\Sigma} (u,v \mid S)|$, and let $\alpha = \beta /2$. Let $M_2$ be the event that Algorithm \ref{Algo:FVS} correctly identifies the $k$-generalized FVS of size $\ell$ is $\cG$.
Then,
\begin{equation}
\P\left( M_2 \right) \geq 1 - \epsilon ,
\end{equation}
and
\begin{align}
&\P\left(\hat{\cE}_{F^c} = \cE_{F^c} \mid F \:
 \text{is a $k$-generalized FVS of size $\ell$} \right) \nonumber \\& \geq 1 - \epsilon ,
\end{align}
for $n = \Omega \left( \frac{\lambda_{\max}^2 + \beta\lambda_{\max}}{\beta^{2}}( (k + \ell) \log p + \log (\epsilon^{-1})) \right)$.
\label{Corr:FVSSample}
\end{corollary}

Just as in the degree bounded and strong separability case, the sample complexity results follows directly from the work on separable graphs \citep{AnaTanHuaWil:12}, so we omit the proof. This algorithm has a computational complexity of $O(p^{k+\ell+4})$, which can be further reduced by using the non-edges in the separability tests to learn other node relations in the graph. Therefore, Algorithm \ref{Algo:FVS} is able to not only identify all the generalized FVSs in a GGM, but also learn their corresponding residual subgraph as well. This is especially useful when the generealized FVS set is small, since most of the edges in the graph can be recovered.


\section{Conclusion}
\label{sect:Conclusion}

In this paper, we introduce two new algorithms for identifying sparse graphs. The first algorithm can identify whether a GGM is $k$-degree bounded or not. In the case of the second algorithm, we introduced the notion of strong separability and showed that the algorithm can determine whether a GGM is strongly $k$-separable or not, with the capability to learn the graph if the graph is indeed strongly $k$-separable. We also establish the concept of a generalized feedback vertex set, and showed that a GGM with a generalized FVS is more general than a strongly separable graph. Finally, we proposed a graph learning algorithm that can identify the generalized FVSs in a GGM, while learning the graph structure outside of the FVS.


\bibliography{latentnodes}

\begin{thebibliography}{18}
\providecommand{\natexlab}[1]{#1}
\providecommand{\url}[1]{\texttt{#1}}
\expandafter\ifx\csname urlstyle\endcsname\relax
  \providecommand{\doi}[1]{doi: #1}\else
  \providecommand{\doi}{doi: \begingroup \urlstyle{rm}\Url}\fi

\bibitem[Anandkumar et~al.(2012)Anandkumar, Tan, Huang, and
  Willsky]{AnaTanHuaWil:12}
A~Anandkumar, V~Y~F Tan, F~Huang, and A~S Willsky.
\newblock High-dimensional gaussian graphical model selection: walk-summability
  and local separation criterion.
\newblock \emph{J. Machine Learning Research}, 13:\penalty0 2293--2337, Aug
  2012.

\bibitem[Belilovskym et~al.(2016)Belilovskym, Varoquaux, and
  Blaschko]{BelVarBla:16}
E~Belilovskym, G~Varoquaux, and M~B Blaschko.
\newblock Testing for differences in gaussian graphical models: Applications to
  brain connectivity.
\newblock In \emph{Advances in Neural Information Processing Systems}, Dec
  2016.

\bibitem[Cicalese and Melani\u{c}(2012)]{CicMil:12}
F~Cicalese and M~Melani\u{c}.
\newblock Graphs of separability at most 2.
\newblock \emph{Discrete Applied Mathematics}, 160\penalty0 (6):\penalty0
  685--696, April 2012.

\bibitem[Durbin et~al.(1999)Durbin, Eddy, Krogh, and
  Mitchison]{DurEddKroMit:99}
R~Durbin, S~R Eddy, A~Krogh, and G~Mitchison.
\newblock \emph{Biological Sequence Analysis: Probabilistic Models of Proteins
  and Nucleic Acids}.
\newblock Cambridge University Press, 1999.

\bibitem[Isard(2003)]{Isa:03}
M~Isard.
\newblock Pampas: real-valued graphical models for computer vision.
\newblock In \emph{IEEE Computer Society Conference on Computer Vision and
  Pattern Recognition}, Jun 2003.

\bibitem[Lauritzen(1996)]{Lau:B96}
S~L Lauritzen.
\newblock \emph{Graphical models}.
\newblock Oxford University Press, New York, 1996.

\bibitem[Liu and Wilsky(2013)]{LiuWil:13}
Y~Liu and A~Wilsky.
\newblock Learning gaussian graphical models with observed or latent fvs.
\newblock In \emph{Advances in Neural Information Processing Systems}, Dec
  2013.

\bibitem[Liu et~al.(2012)Liu, Chandrasekaran, Anandkumar, and
  Willsky]{LiuChaAnaWil:12}
Y~Liu, V~Chandrasekaran, A~Anandkumar, and A~Willsky.
\newblock Feedback message passing for inference in gaussian graphical models.
\newblock \emph{IEEE Trans. Signal Process}, 60\penalty0 (8):\penalty0
  4135--4150, 2012.

\bibitem[Meinshausen and B\"{u}hlmann(2006)]{MeiBuh:06}
N~Meinshausen and P~B\"{u}hlmann.
\newblock High dimensional graphs and variable selection with the lasso.
\newblock \emph{Annals of Statistics}, 34\penalty0 (3):\penalty0 1436--1462,
  2006.

\bibitem[Mohan et~al.(2012)Mohan, Chung, Han, Witten, Lee, and
  Fazel]{MohChuHanWitLeeFaz:12}
K~Mohan, M~J Chung, S~Han, D~Witten, S~Lee, and M~Fazel.
\newblock Structured learning of gaussian graphical models.
\newblock In \emph{Advances in Neural Information Processing Systems}, Dec
  2012.

\bibitem[Pearl(1988)]{Pea:B88}
J.~Pearl.
\newblock \emph{Probabilistic Reasoning in Intelligent Systems}.
\newblock Morgan Kaufmann, 1988.

\bibitem[Ravikumar et~al.(2011)Ravikumar, Wainwright, Raskutti, and
  Yu]{RavWaiRasYu:11}
P~Ravikumar, M~J Wainwright, G~Raskutti, and B~Yu.
\newblock High dimensional covariance estimation by minimizing $\ell$-1
  penalized log-determinant divergence.
\newblock \emph{Electronic Journal in Statistics}, 4:\penalty0 935--980, 2011.

\bibitem[Ren et~al.(2014)Ren, Sun, Zhang, and Zhou]{RenSunZhaZho:14}
Z~Ren, T~Sun, C~Zhang, and H~Zhou.
\newblock Asymptotic normality and optimalities in estimation of large gaussian
  graphical model.
\newblock \emph{Annals of Statistics, to appear}, 2014.

\bibitem[Ryali et~al.(2012)Ryali, Chen, Supekar, and Menon]{Rya:12}
S~Ryali, T~Chen, K~Supekar, and V~Menon.
\newblock Estimation of functional connectivity in fmri data using stability
  selection-based sparse partial correlation with elastic net penalty.
\newblock \emph{Neuroimage}, 59\penalty0 (4):\penalty0 3852--3861, February
  2012.

\bibitem[Soh and Tatikonda(2014)]{SohTat:14}
D~Soh and S~Tatikonda.
\newblock Testing unfaithful gaussian graphical models.
\newblock In \emph{Advances in Neural Information Processing Systems}, Dec
  2014.

\bibitem[Sudderth(2006)]{Sud:06}
E~B Sudderth.
\newblock \emph{Graphical Models for Visual Object Recognition and Tracking}.
\newblock PhD thesis, Massachusetts Institute of Technology. Dept. of
  Electrical Engineering and Computer Science, 2006.

\bibitem[Vazirani(2004)]{Vaz:04}
V~Vazirani.
\newblock \emph{Approximation Algorithms}.
\newblock New York: Springer, 2004.

\bibitem[Whittaker(1990)]{Whi:B90}
J~Whittaker.
\newblock \emph{Graphical Models in Applied Multivariate Statistics}.
\newblock Wiley, 1990.

\end{thebibliography}

\end{document}